\documentclass[twoside,a4paper,12pt,centertags,reqno]{amsart} 
\usepackage{amsmath,amssymb,verbatim,vmargin}
\usepackage{color}
\usepackage{tikz}
\usepackage{hyperref}
\hypersetup{colorlinks,bookmarks=true,linktocpage=true,citecolor=blue,linkcolor=magenta}

\usepackage{eulervm}   

\allowdisplaybreaks 
\usepackage{color}
\usepackage[all]{xy} 

\usepackage{mathtools}
\usepackage{MnSymbol} 

\theoremstyle{plain}
\newtheorem{thm}{Theorem}[section]

\theoremstyle{definition}

\newtheorem{rem}[thm]{Remark}

\usepackage{enumerate}

\newcommand{\bR}{{\mathbb R}}

\newcommand{\cA}{{\mathcal A}}

\newcommand{\cF}{{\mathcal F}}

\def\barint_#1{\mathchoice
            {\mathop{\vrule width 6pt
height 3 pt depth -2.5pt
                    \kern -9.5pt
\intop \kern -4pt}\nolimits_{#1}}%
            {\mathop{\vrule width 5pt height
3 pt depth -2.6pt
                    \kern -6.5pt
\intop \kern -4pt}\nolimits_{#1}}%
            {\mathop{\vrule width 5pt height
3 pt depth -2.6pt
                    \kern -6pt
\intop \kern -4pt}\nolimits_{#1}}%
            {\mathop{\vrule width 5pt height
3 pt depth -2.6pt
          \kern -6pt \intop \kern -4pt}\nolimits_{#1}}}
          
           \def\bariint_#1{\mathchoice
            {\mathop{\vrule width 15pt
height 3 pt depth -2.5pt
                    \kern -15.8pt
\intop \kern -8pt\intop \kern -4pt}\nolimits_{#1}}%
            {\mathop{\vrule width 9pt height
3 pt depth -2.6pt
                    \kern -10.5pt
\intop \kern -8pt\intop \kern -4pt}\nolimits_{#1}}%
            {\mathop{\vrule width 9pt height
3 pt depth -2.6pt
                    \kern -10pt
\intop \kern -8pt\intop \kern -4pt}\nolimits_{#1}}%
            {\mathop{\vrule width 9pt height
3 pt depth -2.6pt
          \kern -8pt \intop \kern -10pt\intop \kern -4pt}
      \nolimits_{  #1}}}

\def\barintlim_#1{\mathchoice
            {\mathop{\vrule width 6pt
height 3 pt depth -2.5pt
                    \kern -8.8pt
\intop \kern -4pt}\limits_{#1}}%
            {\mathop{\vrule width 5pt height
3 pt depth -2.6pt
                    \kern -6.5pt
\intop \kern -4pt}\limits_{#1}}%
            {\mathop{\vrule width 5pt height
3 pt depth -2.6pt
                    \kern -6pt
\intop \kern -4pt}\limits_{#1}}%
            {\mathop{\vrule width 5pt height
3 pt depth -2.6pt
          \kern -6pt \intop \kern -4pt}\limits_{#1}}}
          
           \def\bariintlim_#1{\mathchoice
            {\mathop{\vrule width 15pt
height 3 pt depth -2.5pt
                    \kern -15.8pt
\intop \kern -8pt\intop \kern -4pt}\limits_{#1}}%
            {\mathop{\vrule width 9pt height
3 pt depth -2.6pt
                    \kern -10.5pt
\intop \kern -8pt\intop \kern -4pt}\limits_{#1}}%
            {\mathop{\vrule width 9pt height
3 pt depth -2.6pt
                    \kern -10pt
\intop \kern -8pt\intop \kern -4pt}\limits_{#1}}%
            {\mathop{\vrule width 9pt height
3 pt depth -2.6pt
          \kern -8pt \intop \kern -10pt\intop \kern -4pt}
      \limits_{  #1}}}
          
\renewcommand{\iint}{\int \kern -3pt\int}       

\numberwithin{equation}{section}
\usepackage[symbol]{footmisc}

\makeatletter
\@namedef{subjclassname@2020}{\textup{2020} Mathematics Subject Classification}
\makeatother

\title[Two-point Agmon for graph Schr\"odinger]{A two-point generalisation of the Agmon estimate for Schr\"odinger operators on connected graphs}

\author{Yi C. Huang} 
\address{Yunnan Key Laboratory of Modern Analytical Mathematics and Applications, Yunnan Normal University, Kunming 650500, People's Republic of China}
\address{Department of Mathematical Sciences, Tsinghua University, Beijing 100084, People's Republic of China}
\address{School of Mathematical Sciences, Nanjing Normal University, Nanjing 210023, People's Republic of China}
\email{Yi.Huang.Analysis@gmail.com}
\urladdr{https://orcid.org/0000-0002-1297-7674}

\date{\today} 

\thanks{Research of the author is partially supported by the National NSF grant of China (no. 11801274), 
the Visiting Scholar Program from the Department of Mathematical Sciences of Tsinghua University,
and the Open Project from Yunnan Normal University (no. YNNUMA2403).}

\begin{document}

\begin{abstract}
We provide in this Letter a two-point generalisation of the Agmon estimate for Schr\"odinger operators on graphs recently established by S. Steinerberger.
It reduces to his estimate when the two points belong to different sets separated by the potential and the energy, i.e., the allowed and forbidden regions.
\end{abstract}

\maketitle


\section{Introduction}

Let $G=(V, E)$ be a (finite) connected graph with vertices $V=\{v_1, v_2, \cdots, v_n\}$,
and $(v_i,v_j)\in E$ means there is an edge between $v_i$ and $v_j$ (we simply write $v_i\sim v_j$).
For $v_i\in V$, denote by $\text{deg\,}(v_i)$ the number of vertices connected to $v_i$ by an edge in $E$.

Consider now the following discrete Laplacian $L$ on the graph $G$ (acting on real valued functions $\phi: V\rightarrow\bR$), 
which is formally given by
$$(L\phi)(v)=\text{deg\,}(v)\phi(v)-\sum_{w\sim v}\phi(w), \quad v\in V.$$
Given potential $W: V\rightarrow\bR$ and eigenvalue $E\in\bR$, 
we can further consider the following (spectral) Schr\"odinger equation
\begin{equation} \label{eqn:SchE}
L\phi+W\phi=E\phi.
\end{equation}
The motivation of this letter is to understand the global behaviour of such eigenfunctions $\phi$.
For this purpose we also have to introduce the notion \textit{allowed region} 
$$\cA_V=\{v\in V: W(v)\leq E\},$$ 
the \textit{forbidden region} $\cF_V$ being its complement in $V$.

Given $W$ and $E$, we define the Agmon distance $\rho$ on $(V\times V)\backslash (\cA_V\times \cA_V)$ by 
$$\rho (u,v)=\min\{\rho_1 (u,v),\rho_2 (u,v)\},$$
where
$$\begin{aligned}
\rho_1 (u,v)&=\inf\bigg\{\log\left(1+\frac{(W(v)-E)_+}{\text{deg\,}(v)}\right)+\sum_{i=1}^\ell\log\left(1+\frac{(W(v_i)-E)_+}{\text{deg\,}(v_i)}\right): \\
&\qquad\qquad u=v_1\sim\cdots \sim v_\ell\sim v_{\ell+1}, \,
\frac{(W(v_{\ell+1})-E)_+}{\text{deg\,}(v_{\ell+1})}\geq \frac{(W(v)-E)_+}{\text{deg\,}(v)}\bigg\}
\end{aligned}$$
and 
$$\begin{aligned}
\rho_2 (u,v)&=\inf\bigg\{\log\left(1+\frac{(W(u)-E)_+}{\text{deg\,}(u)}\right)+\sum_{i=1}^\ell\log\left(1+\frac{(W(v_i)-E)_+}{\text{deg\,}(v_i)}\right): \\
&\qquad\qquad v=v_1\sim\cdots \sim v_\ell\sim v_{\ell+1}, \,
\frac{(W(v_{\ell+1})-E)_+}{\text{deg\,}(v_{\ell+1})}\geq \frac{(W(u)-E)_+}{\text{deg\,}(u)}\bigg\},
\end{aligned}$$
and we set
\begin{equation} \label{eqn:Agmon3}
\rho|_{\cA_V\times \cA_V}\equiv0.
\end{equation}
Here as usual, $(a)_+=a$ if $a\geq0$ and $(a)_+=0$ if $a<0$.
For other distances on graphs (combinatorial, resistance, intrinsic, transportation ...), 
see e.g. the monographs by Bakry-Gentil-Ledoux \cite{BakGenLed14},
and by Keller-Lenz-Wojciechowski \cite{KelLenWoj21}.

The following result describes the exponential decay of $\phi$ in terms of $\rho$.

\begin{thm} \label{thm:twopointAgmon}
For $\phi$ a solution of \eqref{eqn:SchE}, and for any $u, v\in V$, we have
\begin{equation} \label{eqn:twopointAgmon}
\min\{|\phi(u)|, |\phi(v)|\}\leq \frac{\|\phi\|_{\ell^\infty(V)}}{e^{\rho (u,v)}}.
\end{equation}
\end{thm}

\begin{rem}
By the maximum principle, we have
$$\|\phi\|_{\ell^\infty(V)}=\|\phi\|_{\ell^\infty(\cA_V)}.$$
Thus for $u\in \cF_V$ and $v\in \cA_V$, \eqref{eqn:twopointAgmon} reduces to Steinerberger's Agmon estimate \cite{Ste23}.
In the continuous setting such estimates can be found in Agmon \cite{Agm14}.
\end{rem}

\begin{rem}
We think the assumption \eqref{eqn:Agmon3} is overlooked in \cite{Ste23}, 
since it is not for sure that two points in $\cA_V$ have to be connected by a path with all nodes in $\cA_V$.
\end{rem}

\section{Proof of Theorem \ref{thm:twopointAgmon}}

The theorem is trivial if both $u$ and $v$ belong to $\cA_V$ as $\rho (u,v)=0$.
So without loss of generality we can assume $u$ is in $\cF_V$ hence $\frac{W(u)-E}{\text{deg\,}(u)}>0$.
The theorem is also trivial if either $\phi(u)$ or $\phi(v)$ is zero.
So we can assume $\min\{|\phi(u)|, |\phi(v)|\}>0$.

First we write down the Schr\"odinger equation on $v_1=u$
$$\text{deg\,}(v_1)\phi(v_1)-\sum_{w\sim v_1}\phi(w)=(E-W(v_1))\phi(v_1), \quad v_1\in V.$$
This can be rewritten as
\begin{equation} \label{eqn:Sch}
\phi(v_1)=\left[1+\frac{W(v_1)-E}{\text{deg\,}(v_1)}\right]^{-1}\frac{1}{\text{deg\,}(v_1)}\sum_{w\sim v_1}\phi(w).
\end{equation}
Note that connectedness is used here.
Taking absolute value on both sides we get
\begin{equation} \label{eqn:Step}
\begin{aligned}
|\phi(v_1)|&\leq\left[1+\frac{W(v_1)-E}{\text{deg\,}(v_1)}\right]^{-1}\frac{1}{\text{deg\,}(v_1)}\sum_{w\sim v_1}|\phi(w)|\\
&\leq\left[1+\frac{W(v_1)-E}{\text{deg\,}(v_1)}\right]^{-1}\max_{w\sim v_1}|\phi(w)|.
\end{aligned}
\end{equation}

Moreover, we deduce from above estimate
$$\max_{w\sim v_1}|\phi(w)|>|\phi(v_1)|\geq \min\{|\phi(u)|, |\phi(v)|\}>0.$$
Hence we can now move from $v_1$ to its neighbour maximising $|\phi(w)|$ (such a vertex is recorded as $v_2$)
and then apply the very same argument again on the vertex $v_2$.
But this argument can be applied iteratively as long as the new vertex still lies in the forbidden region $\cF_V$.
Altogether this results in a path such that
\begin{equation} \label{eqn:stop}
u=v_1\sim\cdots \sim v_m \sim v_{m+1}, \quad \frac{(W(v_{m+1})-E)_+}{\text{deg\,}(v_{m+1})}\geq \frac{(W(v)-E)_+}{\text{deg\,}(v)}.
\end{equation}
Note in particular that $v\in \cA_V\Rightarrow v_{m+1}\in \cA_V$ and $v\in \cF_V\Rightarrow v_{m+1}\in \cF_V$.

Now, collecting the estimates similar to \eqref{eqn:Step} that are produced along the path $u=v_1\sim\cdots \sim v_m \sim v_{m+1}$, 
and if $v\in \cA_V$, we then deduce
$$\begin{aligned}
|\phi(v_1)|&\leq\left(\prod_{i=1}^m\left[1+\frac{W(v_i)-E}{\text{deg\,}(v_i)}\right]^{-1}\right)|\phi(v_{m+1})|\\
&=\left(\prod_{i=1}^m\left[1+\frac{(W(v_i)-E)_+}{\text{deg\,}(v_i)}\right]^{-1}\right)|\phi(v_{m+1})|\\
&\leq \left(\prod_{i=1}^m\left[1+\frac{(W(v_i)-E)_+}{\text{deg\,}(v_i)}\right]^{-1}\right)\|\phi\|_{\ell^\infty(V)}\\
&= \left(\prod_{i=1}^m\left[1+\frac{(W(v_i)-E)_+}{\text{deg\,}(v_i)}\right]^{-1}\right)\left[1+\frac{(W(v)-E)_+}{\text{deg\,}(v)}\right]^{-1}\|\phi\|_{\ell^\infty(V)}.
\end{aligned}$$
If $v\in \cF_V$,
using the equation \eqref{eqn:Sch} for $v_{m+1}$ instead of $v_1$ 
(note that in this case we have in addition $v_{m+1}\in \cF_V$), and the stopping condition from \eqref{eqn:stop}, 
we then deduce
$$\begin{aligned}
|\phi(v_1)|&\leq\left(\prod_{i=1}^m\left[1+\frac{W(v_i)-E}{\text{deg\,}(v_i)}\right]^{-1}\right)|\phi(v_{m+1})|\\
&=\left(\prod_{i=1}^m\left[1+\frac{(W(v_i)-E)_+}{\text{deg\,}(v_i)}\right]^{-1}\right)\\
&\qquad\qquad\times\left[1+\frac{(W(v_{m+1})-E)_+}{\text{deg\,}(v_{m+1})}\right]^{-1}\left|\frac{1}{\text{deg\,}(v_{m+1})}\sum_{w\sim v_{m+1}}\phi(w)\right|\\
&\leq \left(\prod_{i=1}^m\left[1+\frac{(W(v_i)-E)_+}{\text{deg\,}(v_i)}\right]^{-1}\right)\left[1+\frac{(W(v)-E)_+}{\text{deg\,}(v)}\right]^{-1}\|\phi\|_{\ell^\infty(V)}.
\end{aligned}$$
Note that
$$\begin{aligned}
&\left(\prod_{i=1}^m\left[1+\frac{(W(v_i)-E)_+}{\text{deg\,}(v_i)}\right]^{-1}\right)\left[1+\frac{(W(v)-E)_+}{\text{deg\,}(v)}\right]^{-1}\\
&\qquad=\exp\left(-\log\left(1+\frac{(W(v)-E)_+}{\text{deg\,}(v)}\right)-\sum_{i=1}^m\log\left(1+\frac{(W(v_i)-E)_+}{\text{deg\,}(v_i)}\right)\right).
\end{aligned}$$
By the definition of $\rho (u,v)$, we have
$$\exp\left(-\log\left(1+\frac{(W(v)-E)_+}{\text{deg\,}(v)}\right)-\sum_{i=1}^m\log\left(1+\frac{(W(v_i)-E)_+}{\text{deg\,}(v_i)}\right)\right)\leq e^{-\rho (u,v)}.$$
This proves the theorem.

\bigskip

\section*{\textbf{Compliance with ethical standards}}

\bigskip

\textbf{Conflict of interest} The author has no known competing financial interests
or personal relationships that could have appeared to influence this reported work.

\bigskip

\textbf{Availability of data and material} Not applicable.

\bigskip

\bibliographystyle{alpha}

\bibliography{Hua-AgmonTwoPoint} 
 
\end{document}